\documentclass[11 pt]{article}
\usepackage[margin=2.5cm]{geometry}
\usepackage{amsmath,amsthm,amsfonts,amssymb}

\usepackage{graphicx,latexsym}
\usepackage{lscape}
\usepackage[usenames,dvipsnames]{color}
\def\keyw{\par\medskip\noindent {\it \textbf{Keywords}:}\enspace\ignorespaces}
\newtheorem{thm}{Theorem}[section]
 \newtheorem{cor}[thm]{Corollary}
 \newtheorem{lem}[thm]{Lemma}

 \newtheorem{defn}[thm]{Definition}

\newtheorem{rem}[thm]{Remark}

\numberwithin{equation}{section}

\newtheorem{discu}{Discussion:}

\newtheorem{conje}{Conjecture:}

\begin{document}
\date{}

%
%

\title{Embedding complete multi-partite graphs into Cartesian product of paths and cycles}

\author{
R. Sundara\ Rajan $^{a}$
\and
A. Arul Shantrinal $^a$
\and
K. Jagadeesh Kumar $^{a}$
\and
T.M. Rajalaxmi $^{b}$
\and
Jianxi Fan $^{c}$
\and
Weibei Fan $^{c}$
}
\date{}

\maketitle
\vspace{-0.8 cm}
\begin{center}
$^a$ Department of Mathematics, Hindustan Institute of Technology and Science, Chennai, \\ India, 603 103\\
{\tt vprsundar@gmail.com} ~~~~~{\tt shandrinashan@gmail.com}  ~~~~~{\tt jagadeeshgraphs@gmail.com}\\
\medskip

$^b$ Department of Mathematics, SSN College of Engineering, Chennai, India, 603 110\\
{\tt laxmi.raji18@gmail.com}\\
\medskip

$^c$ School of Computer Science and Technology, Soochow University, Suzhou, China, 215006\\
{\tt jxfan@suda.edu.cn} ~~~~~{\tt fanweibei@163.com}\\

\end{center}

\maketitle
\vspace{-0.6 cm}
\begin{abstract}

Graph embedding is a powerful method in parallel computing that maps a guest network $G$ into a host network $H$. The performance of an embedding can be evaluated by certain parameters, such as the dilation, the edge congestion and the wirelength. In this manuscript, we obtain the wirelength (exact and minimum) of embedding complete
multi-partite graphs into Cartesian product of paths and cycles, which include $n$-cube, $n$-dimensional mesh (grid), $n$-dimensional cylinder and $n$-dimensional torus, etc., as the subfamilies.
\end{abstract}

\vspace{-0.1 cm}
\keyw{Embedding, edge congestion, wirelength, complete multi-partite graphs, Cartesian product of graphs}


\vspace{-0.2 cm}
\section{Introduction and Preliminaries}
Given two graphs $G$ (guest) and $H$ (host), an embedding from $G$ to $H$ is an injective mapping $f: V(G) \rightarrow V(H)$ and associating a path $P_f(e)$ in $H$ for each edge $e$ of $G$. We, now define the dilation $dil(G,H)$, the edge congestion $EC(G,H)$ and the wirelength $WL(G,H)$ \cite{BeChHaRoSc98} as follows:
\begin{itemize}
  \item $dil(G,H)=\underset{f:G\rightarrow H}\min ~~\underset{e=xy\in E(G)}\max ~\textrm{dist}_H(f(x),f(y))$
  \item $EC(G,H)=\underset{f:G\rightarrow H}\min ~~\underset{e=xy\in E(H)} \max ~EC_{f}(e)$
  \item $WL(G,H)= \underset{f:G\rightarrow H}\min ~~\underset{e=xy\in E(G)}{\sum } \textrm{dist}_H(f(x),f(y))=\underset{f:G\rightarrow H}\min ~\underset{e=xy\in E(H)}{\sum } ~EC_{f}(e)$
\end{itemize}

where dist$_H(f(x),f(y))$ is a distance (need not be a shortest distance) between $f(x)$ and $f(y)$ in $H$ and $EC_{f}(e)$ denote the number of edges $e'$ of $G$
such that $e=xy$ is in the path $P_{f}(e')$ (need not be a shortest path) between $f(x)$ and $f(y)$ in $H$. Further, $EC_{f}(S)=\underset{e\in S}{\sum }EC_{f}(e)$, where $S \subseteq E(H)$.

\begin{figure}
\centering
\includegraphics[width=10 cm]{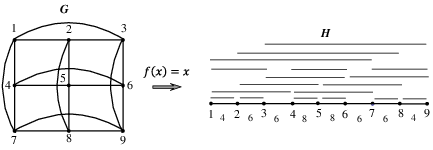}
\caption{An embedding $f$ from torus $G$ into a path $H$ with $dil_f(G,H)=6$, $EC_f(G,H)=8$ and $WL_f(G,H)=48$}
\label{fig1}
\end{figure}

For example, the dilation, the edge congestion and the wirelength of an embedding $f: C_3 \Box C_3 \rightarrow P_9$ is given in Figure \ref{fig1}. It is easy to observe that, the above three parameters are different. But, for any embedding $g$, the sum of the dilations (called dilation sum), the sum of the edge congestion (called edge congestion sum) and the wirelength are all equal. Mathematically, we have the following equality.
$$ \underset{e=xy\in E(G)} \sum \textrm{dist}_H(g(x),g(y))= \underset{e=xy\in E(H)} \sum EC_{g}(e)=WL_g(G,H).$$

In this manuscript, we will use the edge congestion sum to estimate the wirelength.

For a subgraph $M$ of $G$ of order $n$,
\vspace{-0.0 cm}
\begin{itemize}
  \item $I_{G}(M) = \{uv \in E ~|~ u, v \in M\}, ~~I_{G}(k)=\underset{M \subseteq  V(G), ~|M|=k} \max ~|I_{G}(M)|$
  \item $\theta_{G}(M) = \{uv \in E ~|~ u \in M, v  \notin M\}, ~~\theta_{G}(k)=\underset{M \subseteq  V(G), ~|M|=k} \min ~|\theta_{G}(M)|$
 \end{itemize}

The \textit{maximum subgraph problem} (MSP) for a given $k$, $k\in[n]$ is a problem of computing a
subset $M$ of $V(G)$ such that $\left\vert M\right\vert =k$ and $\left\vert I _{G}(M)\right\vert =I _{G}(k)$. Further, the subsets $M$ are
called the \textit{optimal set} \cite{GaJo79, BeDaEl00, Ha04}. Similarly, we define the \textit{minimum cut problem} (MCP) for a given $k$, $k\in[n]$ is a problem of computing a
subset $M$ of $V(G)$ such that $\left\vert M\right\vert =k$ and $\left\vert \theta _{G}(M)\right\vert =\theta _{G}(k)$. For a regular graph, say $r$, we have $2 I_G(k)+\theta_G(k)=rk$, $k \in[n]$ \cite{BeDaEl00}.

The following lemmas are efficient techniques to find the exact wirelength using MSP and MCP.


\begin{lem}{\rm \cite{MiRaPaRa14}}
\label{modifiedcongestionlemma} Let $f:G\rightarrow H$ be an embedding with $|V(G)|=|V(H)|$. Let $S$ be set of all edges (or edge cut) of $H$ such
that $E(H)\setminus S$ has exactly two subgraphs $H_{1}$ and $H_{2}$ and let $G_{i}=[f^{-1}(V(H_{i}))],$ $i=1,2.$ In other words, $G_{i}$ is the induced subgraph on $f^{-1}(V(H_{i}))$ vertices, $i=1,2.$
\textit{Moreover, }$S$ must fulfil the following conditions:
\end{lem}

\vspace{-0.3 cm}
\begin{enumerate}
\item[1.] \textit{For each edge }$ uv\in E(G_{i}), i=1,2$, $P_{f}(uv)$
\textit{has no edges in the set} $S$.

\item[2.] \textit{For each edge }$ uv\in E(G)$ \textit{with} $u$ \textit{in} $V(G_{1})$ \textit{and} $v$ \textit{in} $V(G_{2})$, $P_{f}(uv)$ \textit{has only one edge in the set} $S$.

\item[3.] $V(G_{1})$ \textit{and} $V(G_2)$ \textit{are optimal sets}.
\end{enumerate}

\noindent \textit{Then }$EC_{f}(S)$\textit{\ is minimum over all embeddings $f:G \rightarrow H$ and }\ $EC_{f}(S)=\underset{v\in V(G_1)}{\sum } deg_G(v)-2|E(G_1)|=\underset{v\in V(G_2)}{\sum } deg_G(v)-2|E(G_2)|$, \textit{where $deg_G(v)$ is the degree of a vertex $v$ in $G$}.

\begin{rem}
\rm{For a regular graph $G$, it is easy to check that, $V(G_2)$ is optimal if $V(G_1)$ is optimal and vice-versa \cite{MaRaRaMe09}.}
\end{rem}

\begin{lem} \label{partitionlemma} {\rm \cite{MiRaPaRa14}}
For an embedding $f$ from $G$ into $H$, let $\{P_{1},P_{2},\ldots,P_{t}\}$ be an edge partition of $H$ such
that each $P_{i}$ is an edge cut of $H$ and it satisfies all the conditions of Lemma \ref{modifiedcongestionlemma}. Then%
\begin{equation*}
WL_{f}(G,H)=\overset{t}{\underset{i=1}{\sum }}EC_{f}(P_{i}).
\end{equation*}
\end{lem}

\begin{rem} {\rm \cite {RaRaLiSe18}}
For an embedding $f$ from $G$ into $H$ and $f$ satisfies Lemma \ref{partitionlemma}. Then the wirelength of an embedding from $G$ into $H$ is equal to the wirelength of an embedding from $G$ into $H$ with respect to $f$.
\end{rem}

The multipartite graph is one in all the foremost in style convertible and economical topological structures of interconnection networks. The multipartite has several wonderful options and it's one in all the most effective topological structure of parallel processing and computing systems. In parallel computing, a large process is often decomposed into a collection of little sub processes which will execute in parallel with communications among these sub processes. Due to these communication relations among these sub processes the multipartite graph can be applied for avoiding conflicts in the network as well as multipartite networks helps to identify the errors occurring areas in easy way. A complete $p$-partite graph $G = K_{n_1,\ldots, n_p}$ is a graph that contains $p$ independent sets containing $n_i$, $i\in [p]$, vertices, and all possible edges between vertices from different parts.

The Cartesian product technique is a very powerful technique for create a huger graph from given little graphs and it is very important technique for planning large-scale interconnection networks \cite{Xu01}. Especially, the $n$-dimensional grid (cylinder and torus) structure of interconnection networks offer a really powerful communication pattern to execute a lot of algorithms in many parallel computing systems \cite{Xu01}, which helps to arrange the interconnection network into sequence of sub processors (layers) in uniform distribution manner for transmits the data's in faster way without delay in sending the data packets (messages). Mathematically, we now defined the Cartesian product of graphs as follows:

\begin{defn}{\rm \cite{KhAzAs08}}
Let $G$ and $H$ be two graphs of order $n$ and $m$ respectively. Then the Cartesian product of $G$ and $H$ denoted by $G\Box H$ is the graph with the vertex set $V(G)\times V(H)$, and two vertices $(u,v)$ and $(u',v')$ being adjacent if either $u=u'$ and $vv'\in E(H)$, or $v=v'$ and $uu'\in E(G)$. If $G_1,G_2,\ldots,G_m$ are graphs of order $n_1, n_2, \ldots, n_m$ respectively, then the Cartesian product $G_1\Box G_2\Box \cdots \Box G_m$ is denoted by $\underset{i=1}{\overset{m}\bigotimes}G_i.$
\end{defn}

\begin{rem}
The graph $\underset{i=1}{\overset{n}\bigotimes}G_i$ is said to be an $n$-dimensional grid or torus or cylinder if all $G_i$'s are paths or cycles or any one of the $G_i$ is cycle and the remaining $G_i$'s are paths, respectively.
\end{rem}

The graph embedding problem has been well-studied by many authors with a different networks [1,5,6,7,10--37], and to our knowledge, almost all graphs considered as a host graph is a unique family (for example: path $P_n$, cycle $C_n$, grid $P_n\Box P_m$, cylinder $P_n\Box C_m$, torus $C_n\Box C_m$, hypercube $Q_r$ and so on). In this paper, we overcome this by taking Cartesian product of paths and cycles as a host graph. Moreover, we obtain the wirelength of embedding complete $2^p$-partite graphs $K_{2^{r-p},2^{r-p},\ldots,2^{r-p}}$ into the Cartesian product of graphs $\underset{i=1}{\overset{n}\bigotimes}G_i$, where $G_i$'s are either a path or a cycle, $|V(G_1)|\leq |V(G_2)|\leq \cdots \leq |V(G_n)|$, ~$|V(G_i)|=2^{r_i}$, ~ $|V(\underset{i=1}{\overset{n}\bigotimes}G_i)|=2^{r_1+r_2+\cdots+r_n}=2^r$, ~$r,n\geq 3, ~p\geq 1$ and $1\leq i\leq n$.

\section{Main Result}
To prove the main result, we need the following result and algorithms.

\begin{lem}\label{MSP Complete Partite Graph1}\rm{\cite{RaRaLiSe18}}
If $G$ is a complete $p$-partite graph $K_{r,r,\ldots,r}$ of order $pr$, $p,r\geq 2$, then
\begin{equation*}
I_G(k) =\left\{
\begin{array}{lcl}
\frac{k(k-1)}{2} & ; & k\leq p-1 \\\\
\frac{q^2p(p-1)}{2} & ; & l=qp, 1\leq q \leq r  \\\\
\frac{(q-1)^2p(p-1)}{2}+j(q-1)(p-1)+\frac{j(j-1)}{2} & ; &l=(q-1)p+j, 1\leq j \leq p-1, \\
& & 2\leq q \leq r.
\end{array}
\right.
\end{equation*}
\end{lem}


\noindent \textbf{Guest Graph Algorithm}\\\\
\textbf{Input}: \hspace{0.5cm} $N = 2^{r}$ (Total number of elements)\\
\indent \hspace{1.28 cm}  $p \geq 1$, where $2^{r-p}$ represents the number of elements in the each partite\\\\
\textbf{Output}: \hspace{0.2cm} Labeling of complete $2^{p}$-partite graph $K_{2^{r-p},2^{r-p},\ldots,2^{r-p}}$
\small \begin{enumerate}
  \item Begin the algorithm
  \item The guest graph is generated by the complete $2^{p}$-partite graph
  \item The program contains a function $disp_{-}3nr$ which takes $2^{r-p}, 2^{r-p},\ldots,2^{r-p}$ as partite elements
  \item Get the values $p$ and $N$, where $p\geq1$
  \item $2^{r-p}$ represents the number of elements in a partite
  \item $2^p$ represents the number of partite generated
  \item \textbf{Assign elements in the partite:}
  \item $m=2^{p}$           \hspace{1.5cm} //Determine number of partite
  \item $y = 0$
  \item for $x \leftarrow 0$ to $n$  do
  \item y ++
  \item $Elem_{-}$val = y
  \item for $i \leftarrow 0$ to $p $ do
  \item for $j \leftarrow 0$ to $p$ do
  \item Array[x]+i+j = $Elem_{-}$val
  \item $Elem_{-}$val = $Elem_{-}$val + $N$
  \item Print the partite:
  \item $r = 0 $  \hspace{1.5cm}                                  //Initiating array number
  \item for $x \leftarrow 0$ to $N$ do
  \item for $i \leftarrow 0$ to $p$ do
  \item	 for $z\leftarrow 0$ to $p$ do
  \item for $j \leftarrow 0$ to $p$ do
 \item Print (array [r] + i + j))
  \item Print a tab space
  \item  r++
  \item  Go to new line
  \item   $z = x \% p$
   \item  if $z = 0$
  \item   Print an empty line
\item 	End the algorithm

\end{enumerate}
\normalsize The Python coding and the corresponding implementation of the above algorithm are given in Annexure I.\\


\noindent \textbf{Host Graph Algorithm}\\\\
\textbf{Input}: \hspace{0.5cm} The dimension $n\geq 3$ and the value of $r_1, r_2, \ldots, r_n$\\\\
\textbf{Output}: \hspace{0.0cm} Labeling of Cartesian product of graphs $\underset{i=1}{\overset{n}\bigotimes}G_i$, where $G_i$'s are either a path or a cycle
\small \begin{enumerate}
  \item  The host graph is generated by the Python program hostgraph.py.
  \item  The program contains the main function $disp_{-} n(a_{1},a_{2},\ldots,a_{n})$ which takes parameter $a_{1},a_{2},\ldots,a_{n}$ as a dimension.
  \item  The $a_{1}$ and $a_{2}$ parameter represents the two dimensional base in matrix form.
  \item  The $a_{3}$ parameter takes input which produces $a_{3}$ copies of the base matrix $(a_{1}\times a_{2})$.
  \item  The $a_{4}$ parameter takes input which produces $a_{4}$ copies of the base matrix $(a_{1}\times a_{2}\times a_{3})$.
  \item   Similarly, the $a_{n}$ parameter takes input which produces $a_{n}$ copies of the base matrix \\ $(a_{1}\times a_{2}\times \dots \times a_{n-1})$.
  \item  The number of elements $N$ is given by,
  \begin{center}
  $N =a_{1}\times  a_{2}\times \ldots \times a_{n}$
  \end{center}
  \item The input $a_{1},a_{2},\ldots ,a_{n}$ in the order of $a_{1}\leq a_{2}\leq\ldots \leq a_{n}$.
  \item There will be a $(i-1)$ row-wise rotations for $a_{n}$ copies.
    (i.e) The $i^{th}$ copy of $(a_{1}\times a_{2}\times\ldots \times a_{n-1}\times a_{n})$ say $(a_{1}\times a_{2}\times \ldots \times a_{n-1})^{i}$ has (i - 1) - row wise rotation in the clockwise sense.
  \item  The output is generated by considering $(a_{1}\times a_{2}\times \dots \times a_{n-1})$ base matrix with $a_{n}$ dimension $(a_{1}\times a_{2}\times \dots \times a_{n})$.
\end{enumerate}

\normalsize \noindent The Python coding and the corresponding implementation of the above algorithm are given in Annexure II.

\vspace{0.2 cm}
We, now prove the main result.

\begin{thm}\label{wirelengthalgorithmthm}
Let $G$ be the complete $2^p$-partite graphs $K_{2^{r-p},2^{r-p},\ldots,2^{r-p}}$ and $H$ be the Cartesian product of graphs $\underset{i=1}{\overset{n}\bigotimes}G_i$, where $G_i$'s are either a path or a cycle, $|V(G_1)|\leq |V(G_2)|\leq \cdots \leq |V(G_n)|$, $|V(G_i)|=2^{r_i}$, $|V(\underset{i=1}{\overset{n}\bigotimes}G_i)|=2^{r_1+r_2+\cdots+r_n}=2^r$, $r,n\geq 3, p\geq 1$ and $1\leq i\leq n$. Then the embedding $f$ of $G$ into $H$ given by $f(x)=x$ with minimum (exact) wirelength.
\end{thm}

\noindent \textit{\textbf{Proof.}}~~Label the vertices of $G$ using Guest Graph Algorithm from $1$ to $2^r$. Since the graph $H$ contains an $n$-dimensional grid and label the vertices of $n$-dimensional grid using Host Graph Algorithm from $1$ to $2^r$. For illustration, see Figures 2, 3, 4 and 5. Let us assume that, the label represent each of the vertex, which is allocated by the above algorithms. Let $f:G\rightarrow H$ be an embedding and let $f(v)=v$ for all $v\in V(G)$ and for $uv\in E(G)$, let $P_f(uv)$ be a path (shortest) between $f(u)$ and $f(v)$ in $H$. Now, we have the following 3 cases.

\vspace{-0.2 cm}
\paragraph{Case 1 (All $G_i$'s are paths):}
It is clear that, the graph $H$ becomes an $n$-dimensional grid $P_{2^{r_1}}\Box P_{2^{r_2}}\Box \cdots \Box P_{2^{r_n}}$. For all $i,j$, $1\leq i\leq n$ and $1\leq j\leq 2^{r_i}-1$, let $S_i^j$ be the edge cut of $P_{2^{r_1}}\Box P_{2^{r_2}}\Box \cdots \Box P_{2^{r_n}}$ consisting of the edges between the $j^{th}$ and $j+1^{th}$ copies of \linebreak $P_{2^{r_1}}\Box P_{2^{r_2}}\Box \cdots \Box P_{2^{r_{(i-1)}}}\Box P_{2^{r_{(i+1)}}}\Box \cdots \Box P_{2^{r_n}}$. Then
$\{S_i^j:1\leq i\leq n$ and $1\leq j\leq 2^{r_i}-1\}$ is an edge-partition of $P_{2^{r_1}}\Box P_{2^{r_2}}\Box \cdots \Box P_{2^{r_n}}$.

For all $i,j$, $1\leq i\leq n$ and $1\leq j\leq 2^{r_i}-1$, $E(H)\backslash S_i^j$ has two components $H_{i}^j$ and $\overline{H}_{i}^j$, where
$|V(H_{i}^j)|=(2^{r-r_i})j$ and $|V(\overline{H}_{i}^j)|=2^{r-r_i}(2^{r_i}-j)$. Let $G_{i}^j$ and $\overline{G}_{i}^j$ be the induced subgraph of the inverse images of $V(H_{i}^j)$ and $V(\overline{H}_{i}^j)$ under $f$ respectively. By the guest graph algorithm, $deg_G(v)=2^{r-r_i-p}(2^p-1)j$, for all $v\in V(G_i^j)$ and hence $I_G((2^{r-r_i})j)=2^{2r-2r_i-p}(2^p-1)j^2/2$. By Case 2 of Lemma \ref{MSP Complete Partite Graph1}, $E(G_{i}^j)$ is the maximum subgraph on $|V(G_i^j)|=(2^{r-r_i})j$ vertices in $G$. Thus the edge cut $S_{i}^j$ fulfil all the conditions of Lemma \ref{modifiedcongestionlemma}. Therefore
\begin{eqnarray*}
  EC_{f}(S_{i}^j) &=& 2^{r-p}(2^p-1)(2^{r-r_i})j-2\left(\frac{2^{2r-2r_i-p}(2^p-1)j^2}{2}\right) \\
   &=& 2^{2r-2r_i-p}(2^p-1)(2^{r_i}-j)j
\end{eqnarray*}
is minimum for $1\leq i\leq n$ and $1\leq j\leq 2^{r_i}-1$.

\vspace{0.1 cm}
\noindent Then by Lemma \ref{partitionlemma},
\begin{eqnarray*}
  WL(G, H) &=& \overset{n}{\underset{i=1}{\sum }} ~\overset{2^{r_i}-1}{\underset{j=1}{\sum }}EC_{f}(S_{i}^j)\\
  &=& \overset{n}{\underset{i=1}{\sum }} ~\overset{2^{r_i}-1}{\underset{j=1}{\sum }}2^{2r-2r_i-p}(2^p-1)(2^{r_i}-j)j \\
   &=& \frac{2^{2r-p}(2^p-1)}{6}\left[(2^{r_1}+2^{r_2}+\cdots+2^{r_n})-\left(\frac{1}{2^{r_1}}+\frac{1}{2^{r_2}}+\cdots+\frac{1}{2^{r_n}}\right)\right].
\end{eqnarray*}

\vspace{-0.4 cm}
\paragraph{Case 2 (All $G_i$'s are cycles):}
It is clear that, the graph $H$ becomes an $n$-dimensional torus $C_{2^{r_1}}\Box C_{2^{r_2}}\Box \cdots \Box C_{2^{r_n}}$. For all $i,j$, $1\leq i\leq n$ and $1\leq j\leq 2^{r_i-1}$, let $T_i^j$ be the edge cut of $C_{2^{r_1}}\Box C_{2^{r_2}}\Box \cdots \Box C_{2^{r_n}}$ consisting of the edges between the $(2^{{r_i}-1}-i+j)^{th}$ \& $(2^{{r_i}-1}-i+j+1)^{th}$ and $(2^{r_i}-i+j)^{th}$ \& $(2^{r_i}-i+j+1)^{th}$ copies of $C_{2^{r_1}}\Box C_{2^{r_2}}\Box \cdots \Box C_{2^{r_{(i-1)}}}\Box P_{2^{r_i}-1}\Box C_{2^{r_{(i+1)}}}\Box \cdots \Box C_{2^{r_n}}$. Then
$\{T_i^j:1\leq i\leq n$ and $1\leq j\leq 2^{r_i-1}\}$ is an edge-partition of $C_{2^{r_1}}\Box C_{2^{r_2}}\Box \cdots \Box C_{2^{r_n}}$.

For all $i,j$, $1\leq i\leq n$ and $1\leq j\leq 2^{r_i-1}$, $E(H)\backslash T_i^j$ has two components $H_{i}^j$ and $\overline{H}_{i}^j$, where
$|V(H_{i}^j)|=|V(\overline{H}_{i}^j)|=2^{r-1}$. Let $G_{i}^j$ and $\overline{G}_{i}^j$ be the induced subgraph of the inverse images of $V(H_{i}^j)$ and $V(\overline{H}_{i}^j)$ under $f$ respectively. By the guest graph algorithm, $deg_G(v)=2^{r-p-1}(2^p-1)$, for all $v\in V(G_i^j)$ and hence $I_G(2^{r-1})=2^{2r-2p-2}2^p(2^p-1)/2$. By Case 2 of Lemma \ref{MSP Complete Partite Graph1}, $E(G_{i}^j)$ is the maximum subgraph on $|V(G_i^j)|=2^{r-1}$ vertices in $G$. Thus the edge cut $T_{i}^j$ fulfil all the conditions of Lemma \ref{modifiedcongestionlemma}. Therefore
\begin{eqnarray*}
  EC_{f}(T_{i}^j) &=& 2^{r-p}(2^p-1)2^{r-1}-2\left(\frac{2^{2r-2p-2}2^p(2^p-1)}{2}\right) \\
   &=& 2^{2r-p-2}(2^p-1)
\end{eqnarray*}
is minimum for $1\leq i\leq n$ and $1\leq j\leq 2^{r_i-1}$.

\vspace{0.1 cm}
\noindent Then by Lemma \ref{partitionlemma},
\begin{eqnarray*}
  WL(G, H) &=& \overset{n}{\underset{i=1}{\sum }} ~\overset{2^{r_i-1}}{\underset{j=1}{\sum }}EC_{f}(T_{i}^j)\\
  &=& \overset{n}{\underset{i=1}{\sum }} ~\overset{2^{r_i-1}}{\underset{j=1}{\sum }} 2^{2r-p-2}(2^p-1) \\
   &=& 2^{2r-p-3}(2^p-1)(2^{r_1}+2^{r_1}+\cdots+2^{r_n}).
\end{eqnarray*}

\vspace{-0.5 cm}
\paragraph{Case 3 (Some $G_i$'s are paths and the remaining $G_i$'s are cycles):}

For given $a$ and $b$ of positive integers with $n=a+b$, the host graph $H$ is defined as $G_1\Box G'_1\Box G_2\Box G'_2\Box \cdots \Box G_k\Box G'_q$, where $G_i$'s are path on $2^{r_i}$ vertices or a null graph (it is a graph with $|V|=0$) and $G'_j$'s are cycle on $2^{r_j'}$ vertices or a null graph, $1\leq i\leq k\leq n$, $1\leq j\leq q \leq n$, $k,q \leq n$, $2^{r_1}\leq 2^{r'_1}\leq 2^{r_2}\leq 2^{r'_2}\leq \cdots \leq 2^{r_k}\leq 2^{r'_q}$, $r_1+r_2+\cdots+r_k+r'_1+r'_2+\cdots+r'_q=r$, $r\geq 1$ and $n\geq 3$.

Proceeding along the same lines of Case 1 and Case 2, we describe the edge cuts $S_i^l$ and $T_j^k$, if $G_i$ and $G'_j$ will exists, (we mean $G_i$ and $G'_j$ are not a null graph) respectively, where $1\leq i\leq k\leq n$, $1\leq j\leq q\leq n$, $1\leq l\leq 2^{r_i}-1$ and $1\leq m\leq 2^{r'_j-1}$. Then, it is clear that $\{S_i^l:1\leq i\leq k\leq n, 1\leq l\leq 2^{r_i}-1\}\cup \{T_j^k:1\leq j\leq q\leq n, 1\leq m\leq 2^{r'_j-1}\}$ is an edge-partition of $G_1\Box G'_1\Box G_2\Box G'_2\Box \cdots \Box G_k\Box G'_q$.

\vspace{0.2 cm}
Again by similar arguments in Case 1 and Case 2, we get
\begin{eqnarray*}
  EC_f(S_i^l) &=& 2^{2r-2r_i-p}(2^p-1)(2^{r_i}-l)l
\end{eqnarray*}
is minimum for $1\leq i\leq k$ and $1\leq l\leq 2^{r_i}-1$ and
\begin{eqnarray*}
EC_f(T_j^k) &=& 2^{2r-p-2}(2^p-1)
\end{eqnarray*}
is minimum for $1\leq j\leq q$ and $1\leq m\leq 2^{r'_j-1}$.

\noindent Then by Lemma \ref{partitionlemma},
\begin{eqnarray*}
  WL(G, H) &=& \overset{k}{\underset{\begin{array}{c}
                                       i=1 ~\& \\
                                       G_i ~\textrm{exist}
                                     \end{array}}{\sum }} ~\overset{2^{r_i}-1}{\underset{l=1}{\sum }}~EC_{f}(S_{i}^l) + \overset{q}{\underset{\begin{array}{c}
                                       j=1 ~\& \\
                                       G'_j ~\textrm{exist}
                                     \end{array}}{\sum }} ~\overset{2^{r'_j-1}}{\underset{m=1}{\sum }}~EC_{f}(T_{j}^m)\\
   &=& \overset{k}{\underset{\begin{array}{c}
                                       i=1 ~\& \\
                                       G_i ~\textrm{exist}
                                     \end{array}}{\sum }} ~\overset{2^{r_i}-1}{\underset{l=1}{\sum }}~2^{2r-2r_i-p}(2^p-1)(2^{r_i}-l)l + \overset{q}{\underset{\begin{array}{c}
                                       j=1 ~\& \\
                                       G'_j ~\textrm{exist}
                                     \end{array}}{\sum }} ~\overset{2^{r'_j-1}}{\underset{m=1}{\sum }}~2^{2r-p-2}(2^p-1)\\
  &=& \frac{1}{6}\overset{k}{\underset{\begin{array}{c}
                                       i=1 ~\& \\
                                       G_i ~\textrm{exist}
                                     \end{array}}{\sum }} ~2^{2r-r_i-p}(2^p-1)(2^{r_i}-1) + \overset{q}{\underset{\begin{array}{c}
                                       j=1 ~\& \\
                                       G'_j ~\textrm{exist}
                                     \end{array}}{\sum }} ~2^{2r+r'_j-p-3}(2^p-1). \quad \square
\end{eqnarray*}

\newpage
\begin{cor}
If $G_1$ is a cycle on $2^{r_1}$ vertices and $G_i$ is a path on $2^{r_i}$ vertices, $2\leq i\leq n$, then the host graph $\underset{i=1}{\overset{n}\bigotimes}G_i$ becomes an $n$-dimensional cylinder $C_{2^{r_1}}\Box P_{2^{r_2}}\Box \cdots \Box P_{2^{r_n}}$ and the wirelength of an embedding $f$ from $G$ into $H$ is given by
\begin{eqnarray*}
  WL(G,H) &=& \frac{2^{2r-p}(2^p-1)}{6}\left[(2^{r_2}+2^{r_3}+\cdots+2^{r_n})-\left(\frac{1}{2^{r_2}}+\frac{1}{2^{r_3}}+\cdots+\frac{1}{2^{r_n}}\right)\right]\\
  & & +(2^p-1)2^{2r+r'_1-p-3}. \quad \square
\end{eqnarray*}
\end{cor}

\vspace{-0.4 cm}
\section{Conclusion and Future Work}
In this manuscript, we found the wirelength (exact and minimum) of an embedding complete multi-partite graphs into Cartesian product of paths and cycles. Computing the dilation and the edge congestion of embedding complete multi-partite graphs into Cartesian product and other product of graphs are under investigation.

\vspace{6 pt} \noindent
{\bf Acknowledgments}\
\vspace{0.4 cm}

The work of R. Sundara Rajan was partially supported by Project no. ECR/2016/1993, Science and Engineering Research Board (SERB), Department of Science and Technology (DST), Government of India. Further, we thank Prof. Gregory Gutin and Prof. Stefanie Gerke, Royal Holloway, University of London, TW20 0EX, UK; Prof. Indra Rajasingh, School of Advanced Sciences, Vellore Institute of Technology, Chennai, India and Dr. N. Parthiban, Department of Computer Science and Engineering, SRM Institute of Science and Technology, Chennai, India, for their fruitful suggestions.

\newpage
\begin{center}
\textbf{Annexure I}
\end{center}

\noindent \textbf{Python program for labeling of the guest graph}
\begin{flushleft}
def printline(num, boxes, boxesinrow):\\
\hspace{0.8cm}    i = 0;\\
\hspace{0.8cm}    string = ' '\\
\hspace{0.8cm}    for i in range(0, boxesinrow):\\
\hspace{1cm}        for j in range(0, 4):\\
\hspace{1.4cm}            string = string + "$\{:<3d\}$ ".format(num + boxes * j) + '  '\\
\hspace{1cm}        string = string + '       '\\
\hspace{1cm}        num = num + 1\\
\hspace{1cm}    print(string)\\

def printbox(num, boxes, $elements_{-}$ $per_{-}$box, $boxes_{-}$ $in_{-}$row):\\
\hspace{0.8cm}    value=num\\
\hspace{0.8cm}    temp = $elements_{-}$ $per_{-}$ box // 4\\
\hspace{0.8cm}    while temp $>$ 0:\\
\hspace{1cm}        printline(num, boxes, $boxes_{-}$ $in_{-}$ row)\\
\hspace{1cm}        temp = temp - 1\\
\hspace{1cm}        num = num + boxes * 4\\
\hspace{0.8cm}    print('\textbackslash n') \\
\hspace{0.8cm}    return (value+$boxes_{-}$ $in_{-}$ row)\\

def printpattern(numl):\\
\hspace{0.8cm}    boxes = len(numl)\\
\hspace{0.8cm}    $elements_{-}$ $Per_{-}$ box = numl[0]\\
\hspace{0.8cm}    num = 1\\
\hspace{0.8cm}    temp = boxes\\
\hspace{0.8cm}    while $temp > 0$:\\
\hspace{1cm}        if $temp >= 4$:\\
\hspace{1.4cm}            num = printbox(num, boxes, $elements_{-}$ $Per_{-}$ box, 4)\\
\hspace{1cm}        else:\\
\hspace{1.4cm}            num = printbox(num, boxes, $elements_{-}$ $Per_{-}$ box, temp)\\
\hspace{1cm}        temp = temp - 4\\

printpattern([32,32,\ldots,32])\\
\end{flushleft}

\newpage
\noindent \textbf{Implementation of the above Python program}\\\\
\noindent \textbf{Output 1:}

\vspace{-0.2 cm}
$$\includegraphics[width=16 cm]{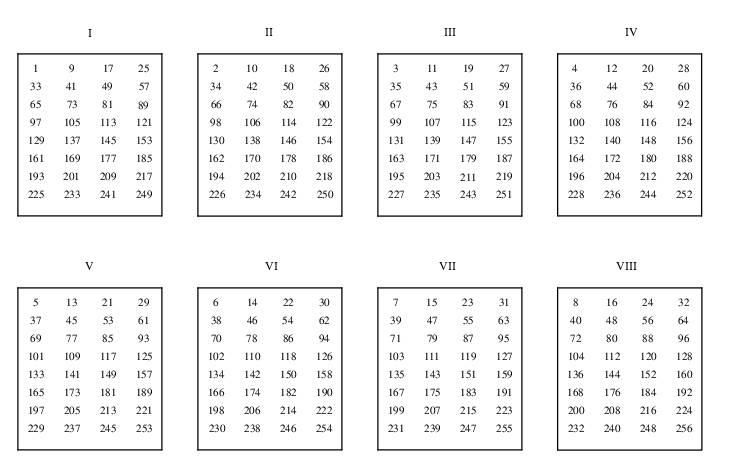}
\label{fig2}$$
\vspace{-0.5 cm}
\begin{center}
Figure 2: Complete 8-partite graph $K_{{32},32,\ldots,{32}}$
\end{center}

\newpage
\noindent \textbf{Output 2:}
\vspace{-0.2 cm}
$$\includegraphics[width=15 cm]{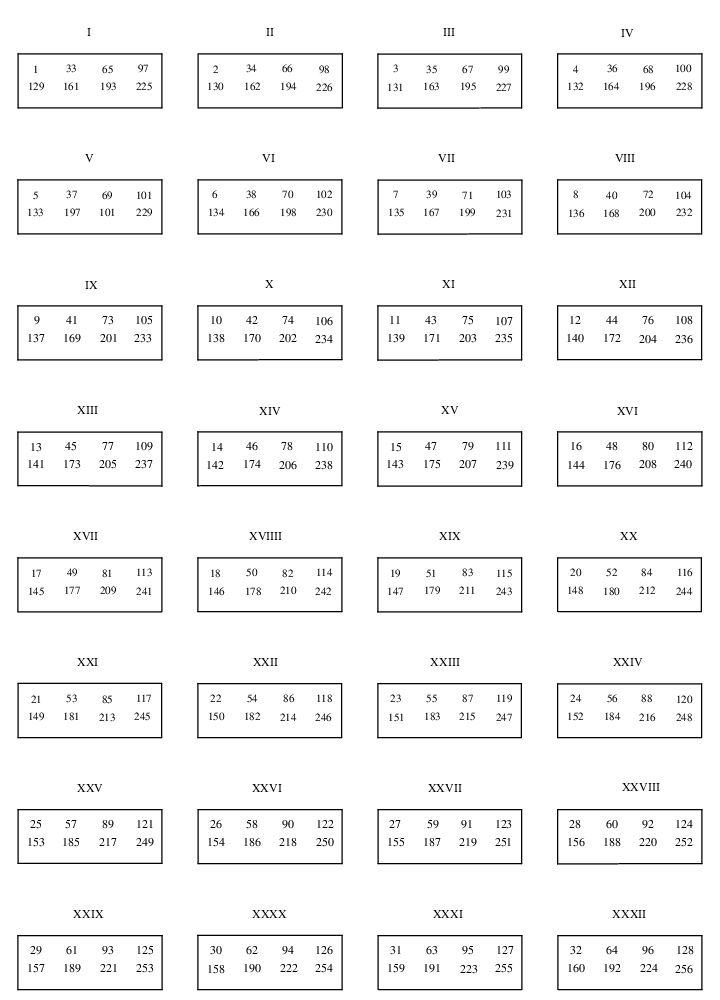}
\label{fig3}$$
\vspace{-0.5 cm}
\begin{center}
Figure 3: Complete 32-partite graph $K_{{8},8,\ldots,{8}}$
\end{center}

\newpage
\newpage
\begin{center}
\textbf{Annexure II}
\end{center}
\vspace{-0.1 cm}
\small
\noindent \textbf{Python program for labeling of the host graph}
\vspace{-0.2 cm}
\begin{flushleft}
n=0

def $disp_{-}$3nr(numl, n, irange):\\
\hspace{0.5cm}    for i in irange:\\
\hspace{1cm}         string = str(i+n)\\
\hspace{1cm}        for j in range(1, numl[1]):\\
\hspace{1.4cm}            string=string+' '+str(i+n+numl[0]*j)\\
\hspace{1cm}        for k in range(1, numl[2]):\\
\hspace{1.4cm}             string=string+'     '\\
\hspace{1cm}            for j in range(0, numl[1]):\\
\hspace{1.7cm}                string = string + '    ' + str(i+ n + numl[0] * j+ k*numl[0]*numl[1])\\
\hspace{1cm}        print(string)

def $disp_{-}$ n(numl):\\
\hspace{0.5cm}    base=numl[0]*numl[1]*numl[2]\\
\hspace{0.5cm}    $para_{-}$ num=len(numl)\\
\hspace{0.5cm}    order=numl[3:]\\
\hspace{0.5cm}    global n\\
\hspace{0.5cm}    n = -base\\
\hspace{0.5cm}    tnum1=numl[:3]\\
\hspace{0.5cm}    loopri(order, tnum1, base)\\

def rotate(irange):\\
\hspace{0.5cm}    temp=irange[0]\\
\hspace{0.5cm}    for i in range(0, len(irange)-1):\\
\hspace{1cm}        irange[i] = irange[i+1]\\
\hspace{0.5cm}    irange[len(irange)-1]=temp\\
\hspace{0.5cm}    return(irange)\\

def loopri(order, tnuml, base):\\
\hspace{0.5cm}    irange = list(range(1, tnuml[0]+1))\\
\hspace{0.5cm}    if len(order) == 1:\\
\hspace{1cm}         for i in range(0, order[0]):\\
\hspace{1.4cm}             global n\\
\hspace{1.4cm}            print(i+1)\\
\hspace{1.4cm}            n = n + base\\
\hspace{1.4cm}            $disp_{-}$ 3nr(tnuml, n, irange)\\
\hspace{1.4cm}            irange=rotate(irange)\\
\hspace{1.4cm}            print('   ')\\
\hspace{0.5cm}    else:\\
\hspace{1cm}        for i in range(0, order[-1]):\\
\hspace{1.4cm}            print(i)\\
\hspace{1.4cm}            loopr(order[:-1], tnuml, base, irange)\\
\hspace{1.4cm}            irange=rotate(irange)\\

def loopr(order, tnuml, base, irange):\\
\hspace{0.5cm}    if len(order) == 1:\\
\hspace{1cm}        for i in range(0, order[0]):\\
\hspace{1.4cm}            global n\\
\hspace{1.4cm}            n = n + base\\
\hspace{1.4cm}            $disp_{-}$ 3nr(tnuml, n, irange)\\
\hspace{1.4cm}            print('     ')\\
\hspace{0.5cm}    else:\\
\hspace{1cm}         for i in range(0, order[-1]):\\
\hspace{1.4cm}            print(i)\\
\hspace{1.4cm}            loopr(order[:-1], tnuml, base, irange)\\

$disp_{-}$ n([4,8,16])\\
\end{flushleft}

\noindent \textbf{Implementation of the above Python program}\\\\
\noindent \textbf{Output 1:}

\vspace{-0.2 cm}
$$\includegraphics[width=9 cm, height=20.8 cm]{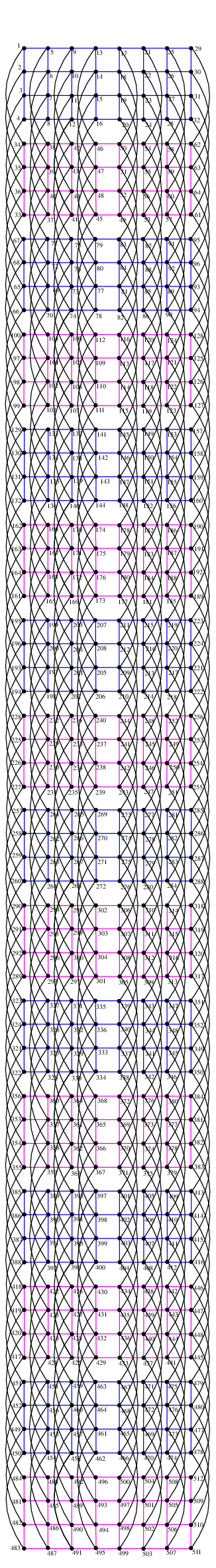}
\label{fig4}$$
\vspace{-0.5 cm}
\begin{center}
Figure 4: 3-dimensional grid $P_4\Box P_8\Box P_{16}$
\end{center}

\noindent \textbf{Output 2:}
\vspace{-0.2 cm}
$$\includegraphics[width=10 cm, height=21.8 cm]{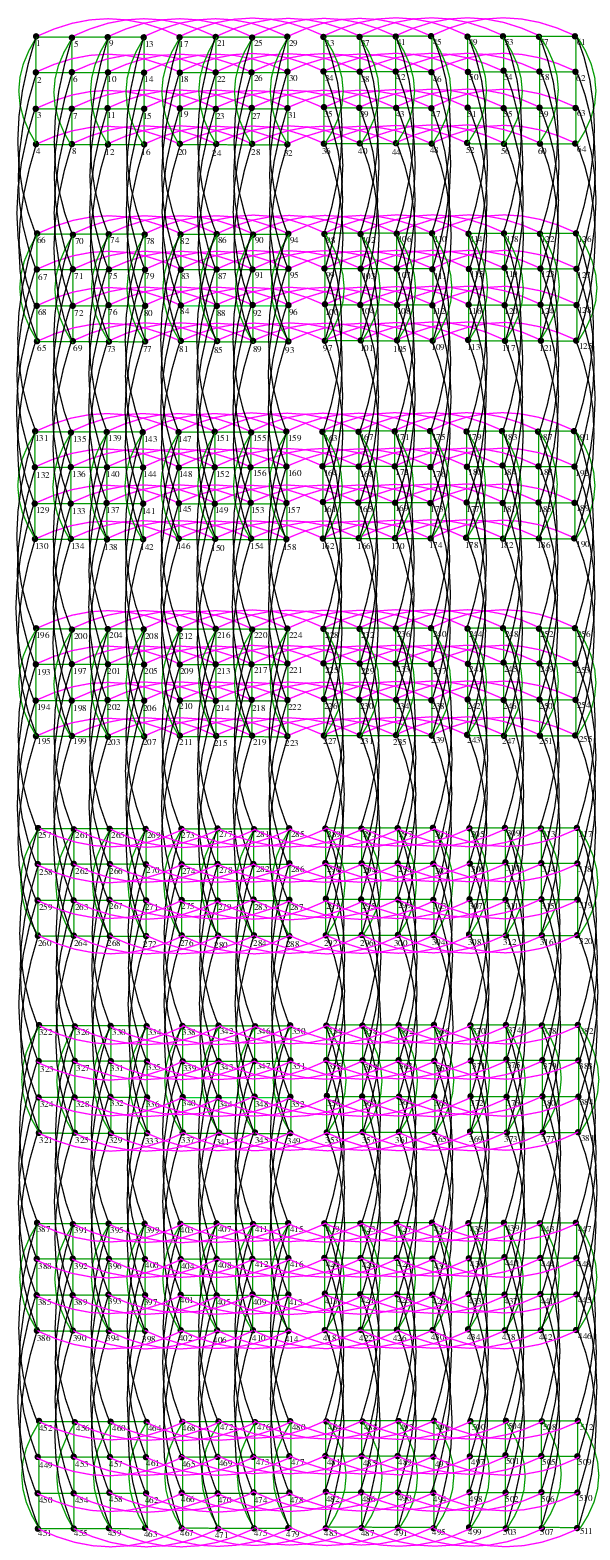}
\label{fig5}$$
\vspace{-0.5 cm}
\begin{center}
Figure 5: 4-dimensional cylinder $C_{4}\Box P_{4}\Box P_{4}\Box P_{8}$
\end{center}

\end{document}